\documentclass[a4paper]{book}
\usepackage{amsmath,latexsym,amssymb,kcp2}

\newcommand{\na}[1]{\mathit{#1}}    
\newcommand{\fn}[1]{\mathrm{#1}}    

\newcommand{\limplies}{\rightarrow}
\newcommand{\liff}{\leftrightarrow}
\newcommand{\ex}[1]{\exists #1 \;} 
\newcommand{\fa}[1]{\forall #1 \;} 
\newcommand{\ph}{\varphi}

\newcommand{\la}{(}
\newcommand{\ra}{)}

\newcommand{\N}{\mathsf{N}}

\newcommand{\PR}{\fn{PR}^\omega}
\newcommand{\lam}[1]{\lambda #1 \;} 

\newcommand{\inot}{\mathord{\sim}}

\newcommand{\isleft}[1]{\mathrm{isleft}(#1)}
\newcommand{\isright}[1]{\mathrm{isright}(#1)}

\newcommand{\eltl}[1]{\mathrm{left}(#1)}
\newcommand{\eltr}[1]{\mathrm{right}(#1)}
\newcommand{\inl}{\mathrm{inl}}
\newcommand{\inr}{\mathrm{inr}}

\newcommand{\nnf}{\fn{nnf}}

\newcommand{\awk}{\mathit{awk}}
\newcommand{\realizes}{\mathop{\; \mathit{realizes} \;}}

\collection

\mainmatter

\begin{document}

\paper{The Computational Content of Classical Arithmetic\thanks{Dedicated to Grigori Mints in honor of his seventieth birthday.}}{Jeremy Avigad}

\begin{abstract}
  Almost from the inception of Hilbert's program, foundational and
  structural efforts in proof theory have been directed towards the
  goal of clarifying the computational content of modern mathematical
  methods. This essay surveys various methods of extracting
  computational information from proofs in classical first-order
  arithmetic, and reflects on some of the relationships between
  them. Variants of the G\"odel-Gentzen double-negation translation,
  some not so well known, serve to provide canonical and efficient
  computational interpretations of that theory.
\end{abstract}

\section{Introduction}

Hilbert's program was launched, in 1922, with the specific goal of
demonstrating the consistency of modern, set-theoretic methods, using
only finitary means. But the program can be viewed more broadly as a
response to the radical methodological changes that had been
introduced to mathematics in the late nineteenth century. Central to
these changes was a shift in mathematical thought whereby the goal of
mathematics was no longer viewed as that of developing powerful
methods of calculation, but, rather, that of characterizing abstract,
possibly infinite, mathematical structures, often in ways that could
not easily be reconciled with a computational understanding.

Grisha's work over the years has touched on almost every aspect of
proof theory, both of the reductive (foundational) and structural
sort, involving a wide range of logical frameworks. But much of his
work addresses the core proof-theoretic concern just raised, and has
served to provide us with a deep and satisfying understanding of the
computational content of nonconstructive, infinitary reasoning. Such
work includes his characterization of the provably total computable
functions of $\na{I\Sigma_1}$ as exactly the primitive recursive
functions (\cite{mints:71,mints:76}); his method of continuous cut
elimination, which provides a finitary interpretation of infinitary
cut-elimination methods (\cite{mints:92b,buchholz:91}); and his work
on the epsilon substitution method (for example,
\cite{mints:96,mints:tupailo:99}).

Grisha has also been a friend and mentor to me throughout my career.
The characterization of the provably total computable functions of
$\na{I\Sigma_1}$ just mentioned was, in fact, also discovered by
Charles Parsons and Gaisi Takeuti, all independently. I shudder to
recall that at a meeting at Oberwolfach in 1998, when I was just
two-and-a-half years out of graduate school, I referred to the result
as ``Parsons' theorem'' in a talk before an audience that, unfortunately,
included only the other two. Grisha asked the first question after the
talk was over, and nothing in his manner or tone even hinted that I
had made a \emph{faux pas} (it didn't even occur to me until much
later). In fact, I still vividly remember his encouraging and
insightful comments, then and in later discussion. (For the record,
Gaisi was equally gracious and supportive.)

In this essay, I will discuss methods of interpreting classical
first-order arithmetic, often called \emph{Peano arithmetic}
($\na{PA}$), in computational terms. Although the study of $\na{PA}$
was central to Hilbert's program, it may initially seem like a toy
theory, or an artificially simple case study. After all, mathematics
deals with much more than the natural numbers, and there is a lot more
to mathematical argumentation than the principle of induction. But
experience has shown that the simplicity of the theory is deceptive:
via direct interpretation or more elaborate forms of proof-theoretic
reduction, vast portions of mathematical reasoning can be understood
in terms of $\na{PA}$ \cite{avigad:03b,feferman:87b,simpson:99}.

Here, I will be concerned with the $\Pi_2$, or ``computational,''
consequences of $\na{PA}$. Suppose $\na{PA}$ proves $\fa x \ex y
R(x,y)$, where $x$ and $y$ range over the natural numbers and $R(x,y)$
is a decidable (say, primitive recursive) relation. We would like to
understand how and to what extent such a proof provides an
\emph{algorithm} for producing such a $y$ from a given $x$, one that
is more informative than blind search. There are four methods that are
commonly used to extract such an algorithm:
\begin{enumerate}
\item G\"odel's \emph{Dialectica} interpretation
  \cite{goedel:58,avigad:feferman:98}, in
  conjunction with a double-negation interpretation that interprets
  $\na{PA}$ in its intuitionistic counterpart, Heyting arithmetic
  ($\na{HA}$)
\item realizability
  \cite{kleene:45,kleene:52,kreisel:59,troelstra:98}, again in
  conjunction with a double-negation translation, and either the
  Friedman $A$-translation \cite{friedman:78} (often also attributed
  to Dragalin and Leivant, independently) or a method due to Coquand
  and Hofmann (\cite{coquand:hofmann:99,avigad:00b}) to ``repair''
  translated $\Pi_2$ assertions
\item cut elimination (\cite{gentzen:36}; see, for example,
  \cite{schwichtenberg:77}) 
\item the epsilon substitution method
  (\cite{hilbert:bernays:34,avigad:zach:02})
\end{enumerate}
These four approaches really come in two pairs: the Dialectica
interpretation and realizability have much in common, and, indeed,
Paulo Oliva \cite{oliva:06} has recently shown that one can
interpolate a range of methods between the two; and, similarly, cut
elimination and the epsilon substitution method have a lot in common,
as work by Grisha (e.g.~\cite{mints:08}) shows. That is not to say
that there aren't significant differences between the methods in each
pairing, but the differences between the two pairs are much more
pronounced.

For one thing, they produce two distinct sorts of ``algorithms.'' The
Dialectica interpretation, and Kreisel's ``modified'' version of
Kleene's realizability, extract terms in G\"odel's calculus of
primitive recursive functionals of finite type, denoted $\PR$ in
Section~\ref{preliminaries:section} below. In contrast, cut
elimination and the epsilon substitution method provide iterative
procedures, whose termination can be proved by ordinal analysis.
Specifically, one assigns (a notation for) an ordinal less than
$\varepsilon_0$ to each stage of the computation in such a way that
the ordinals decrease as the computation proceeds. Terms in $\PR$ and
$\mathord{\prec}\varepsilon_0$-recursive algorithms both have
computational meaning, and there are various ways to ``see'' that the
computations terminate; but, of course, any means of proving
termination formally for all such terms and algorithms has to go
beyond the means of $\na{PA}$.

Second, as indicated above, the first two methods involve an
intermediate translation to $\na{HA}$, while the second two do not. It
is true that the Dialectica interpretation and realizability can be
applied to classical calculi directly (see \cite{shoenfield:01} for
the Dialectica interpretation, and, for example,
\cite{avigad:00,murthy:91} for realizability); but I know of no such
interpretation that cannot be understood in terms of a passage through
intuitionistic arithmetic
\cite{avigad:00,avigad:unp:k,streicher:kohlenbach:07}. In contrast,
cut elimination and the epsilon substitution method apply to classical
logic directly. That is not to deny that one can apply cut elimination
methods to intuitionistic logic (see, for example,
\cite{troelstra:schwichtenberg:00}); but the arguments tend to be
easier and more natural in the classical setting.

Finally, there is the issue of canonicity. Algorithms extracted from
proofs in intuitionistic arithmetic tend to produce canonical
witnesses to $\Pi_2$ assertions; work by Grisha
\cite{mints:92c,mints:92d} shows, for example, that algorithms
extracted by various methods yield the same results. In contrast,
different ways of extracting witnesses from classical proofs yield
different results, conveying the impression that there is something
``nondeterministic'' about classical logic. (There is a very nice
discussion of this in \cite{urban:00,urban:bierman:01}. See also the
discussion in Section~\ref{conclusions:section} below.) Insofar as one
has a natural translation from classical arithmetic to intuitionistic
arithmetic, some of the canonicity of the associated computation is
transferred to the former theory.

In this essay, I will discuss realizability and the Dialectica
interpretation, as they apply to classical arithmetic, via
translations to intuitionistic arithmetic. After reviewing some
preliminaries in Section~\ref{preliminaries:section}, I will discuss
variations of the double-negation interpretation in
Section~\ref{double:negation:section}. One, in particular, is very
efficient when it comes to introducing negations; in
Section~\ref{realizability:section}, I will show that, when combined
with realizability or the Dialectica interpretation, this yields
computational interpretations of classical arithmetic that are
efficient in their use of higher types. In
Section~\ref{herbrand:section}, I will consider another curious
double-negation interpretation, and diagnose an unfortunate aspect of
its behavior.

This work has been partially supported by NSF grant DMS-0700174 and a
grant from the John Templeton Foundation. I am grateful to Philipp Gerhardy, Thomas Streicher, and an anonymous
referee for helpful comments and corrections. 

\section{Preliminaries}
\label{preliminaries:section}

Somewhat imprecisely, one can think of intuitionistic logic as
classical logic without the law of the excluded middle; and one can
think of minimal logic as intuitionistic logic without the rule
\emph{ex falso sequitur quodlibet}, that is, from $\bot$ conclude
anything. Computational interpretations of classical logic often pass
through minimal logic, which has the nicest computational
interpretation. (One can interpret intuitionistic logic in minimal
logic by replacing every atomic formula $A$ by $A \lor \bot$, so the
difference between these two is small.)

To have a uniform basic to compare the different logics, it is useful
to take the first-order logical symbols to be $\forall$, $\exists$,
$\land$, $\lor$, $\limplies$, and $\bot$, with $\lnot \ph$ defined to
be $\ph \limplies \bot$. However, when it comes to classical logic, it
is often natural to restrict one's attention to formulas in
\emph{negation-normal form}, where formulas are built up from atomic
and negated atomic formulas using $\land$, $\lor$, $\forall$, and
$\exists$. A negation operator, $\inot \ph$, can be defined for such
formulas; $\inot \ph$ is what you get if, in $\ph$, you exchange
$\land$ with $\lor$, $\forall$ with $\exists$, and atomic formulas
with their negations. Note that $\inot \inot \ph$ is just
$\ph$. Classically, every formula $\ph$ has a negation-normal form
equivalent, $\ph^\nnf$, obtained by defining $(\theta \limplies
\eta)^\nnf$ to be $\inot \theta^\nnf \lor \eta^\nnf$, and treating the
other connectives in the obvious way. This has the slightly awkward
consequence that $(\lnot \ph)^\nnf$ translates to $\inot \ph^\nnf \lor
\bot$, but simplifying $\theta \lor \bot$ to $\theta$ and $\theta
\land \bot$ to $\bot$ easily remedies this.

There are a number of reasons why negation-normal form is so natural
for classical logic. First of all, it is easy to keep track of
polarities: if $\ph$ is in negation-normal form, then every subformula
is a positive subformula, except for, perhaps, atomic formulas; an
atomic formula $A$ occurs positively in $\ph$ if it occurs un-negated,
and negatively if it occurs with a negation sign before it. Second,
the representation accords well with practice: any classically-minded
mathematician would not hesitate to prove ``if $\ph$ then $\psi$'' by
assuming $\lnot \psi$ and deriving $\lnot \ph$, or by assuming $\ph$
and $\lnot \psi$ and deriving a contradiction; so it is convenient
that $\ph \limplies \psi$, $\lnot \psi \limplies \lnot \ph$, and
$\lnot (\ph \land \lnot \psi)$ have the same negation-normal form
representation. Finally, proof systems for formulas in negation-normal
form tend to be remarkably simple (see, for example,
\cite{schwichtenberg:77,troelstra:schwichtenberg:00}).

It was G\"odel \cite{goedel:58} who first showed that the provably
total computable functions of arithmetic can be characterized in terms
of the primitive recursive functionals of finite type (see
\cite{avigad:feferman:98,hindley:seldin:86}). The set of
finite types can be defined to be the smallest set containing the
symbol $\N$, and closed under an operation which takes types $\sigma$
and $\tau$ to a new type $\sigma \to \tau$. In the intended (``full'')
interpretation, $\N$ denotes the set of natural numbers, and $\sigma
\to \tau$ denotes the set of all functions from $\sigma$ to $\tau$.
A set of terms, $\PR$, is defined inductively as follows:
\begin{enumerate}
\item For each type $\sigma$, there is a stock of variables $x, y,
  z,\ldots$ of type $\sigma$.
\item $0$ is a term of type $\N$.
\item $S$ (successor) is a term of type $\N \to \N$.
\item if $s$ is a term of type $\tau \rightarrow \sigma$ and $t$ is a
  term of type $\tau$, then $s(t)$ is a term of type $\sigma$.
\item if $s$ is a term of type $\sigma$ and $x$ is a variable of type
  $\tau$, then $\lam x s$ is a term of type $\tau \rightarrow
  \sigma$.
\item If $s$ is a term of type $\sigma$, and $t$ is a term of type
  $\N \to (\sigma \to \sigma)$, then $R_{st}$ is a term of type
  $\N \to \sigma$.
\end{enumerate}
Intuitively, $s(t)$ denotes the result of applying $s$ to $t$, $\lam x
s$ denotes the function which takes any value of $x$ to $s$, and
$R_{st}$ denotes the function defined from $s$ and $t$ by primitive
recursion, with $R_{st}(0) = s$ and $R_{st}(S(x)) = t(x,R_{st}(x))$
for every $x$.
In this last equation, I have adopted the convention of writing
$t(r,s)$ instead of $(t(r))(s)$.

It will be convenient below to augment the finite types with products
$\sigma \times \tau$, associated pairing operations $\la \cdot, \cdot
\ra$, and projections $(\cdot)_0$ and $(\cdot)_1$.  Product types can
be eliminated in the usual way by currying and replacing terms $t$
with sequences of terms $t_i$. It will also be convenient to have
disjoint union types $\sigma + \tau$, an element of which is either an
element of $\sigma$ or an element of $\tau$, tagged to indicate which
is the case. That is, for each such type we have insertion operations,
$\inl$ and $\inr$, which convert elements of type $\sigma$ and $\tau$
respectively to an element of type $\sigma + \tau$; predicates
$\isleft{a}$ and $\isright{a}$, which indicate whether $a$ is tagged
to be of type $\sigma$ or $\tau$; and functions $\eltl{a}$ and
$\eltr{a}$, which interpret $a$ as an element of type $\sigma$ and
$\tau$, respectively. References to such sum types can be eliminated
by taking $\sigma + \tau$ to be $\N \times \sigma \times \tau$,
defining $\inl(a) = \la 0, a, 0^\tau \ra$, defining $\inr(a) = \la 1,
0^\sigma, a \ra$, where $0^\sigma$ and $0^\tau$ are constant zero
functionals of type $\sigma$, $\tau$ respectively, and so on.


In the next section, I will describe various double-negation
interpretations that serve to reduce classical arithmetic, $\na{PA}$,
to intuitionistic arithmetic, $\na{HA}$ --- in fact, to $\na{HA}$
taken over \emph{minimal logic}. These show that if $\na{PA}$ proves a
$\Pi_2$ formula $\fa x \ex y R(x,y)$, then $\na{HA}$ proves $\fa x
\lnot \lnot \ex y R(x,y)$; in fact, a variant $\na{HA'}$ of $\na{HA}$
based on minimal logic suffices. This reduces the problem to
extracting computational information from the latter proof.

One method of doing so involves using Kreisel's notion of modified
realizability, combined with the Friedman A-translation. One can
extend $\na{HA}$ to a higher-type version, $\na{HA^\omega}$, which has
variables ranging over arbitrary types, and terms of all the primitive
recursive functionals. Fix any primitive recursive relation $A(y)$;
then to each formula $\ph(\bar x)$ in the language of arithmetic, one
inductively assigns a formula ``$a \realizes \ph(\bar x)$,'' as
follows.
\begin{alignat*}{2}
  & a \realizes \bot & \; & \equiv \; A(a) \\
  & a \realizes \theta & \; & \equiv \; \theta, \; \mbox{if $\theta$ is atomic} \\
  & a \realizes \ph \land \psi & \; & \equiv \; ((a)_0 \realizes \ph)
  \land
  ((a)_1 \realizes \ph) \\
  & a \realizes \ph \lor \psi & \; & \equiv \; ((\isleft{a} \land
  \eltl{a} \realizes
  \ph) \mathop{\lor} \\
  & & \; & \quad \quad \quad (\isright{a} \land \eltr{a} \realizes \psi)) \\
  & a \realizes \ph \limplies \psi & \; & \equiv \; \fa b (b
  \realizes \ph \limplies a(b) \realizes \psi) \\
  & a \realizes \fa x \ph(x) & \; & \equiv \; \fa x (a(x) \realizes \ph(x)) \\
  & a \realizes \ex x \ph(x) & \; & \equiv \; (a)_1 \realizes
  \ph((a)_0)
\end{alignat*}
Now, suppose classical arithmetic proves $\fa x \ex y R(x,y)$, for
some primitive recursive relation $R$. Then, using a double-negation
translation, $\na{HA'}$ proves $\fa x \lnot \lnot \ex y R(x,y)$, and
hence it proves $\lnot \lnot \ex y R(c,y)$ for a fresh constant
$c$. Fix $A(y)$ in the realizability relation above to be the formula
$R(c,y)$. Inductively, one can then extract from the proof of term $t$
of $\PR$ such that $\na{HA^\omega}$ proves that $t$ realizes $\lnot
\lnot \ex y R(x,y)$. Now notice that the identity function, $\na{id}$,
realizes $\lnot \ex y R(c,y)$, since a realizer to $\ex y R(c,y)$ is
simply a value of $a$ satisfying $R(c,a)$. Thus if $a$ realizes $\ex y
R(c,y)$, then $a(\na{id})$ satisfies $R(c,a(\na{id}))$. Viewing $a$,
now, as a function of $c$, yields the following conclusions:
\begin{theorem}
\label{main:witnessing:thm}
If classical arithmetic proves $\fa x \lnot \lnot \ex y R(x,y)$, there is a
term $F$ of $\PR$ of type $\N$ to $\N$ such that
$\na{HA^\omega}$ proves $\fa x R(x,F(x))$.
\end{theorem}
\noindent See \cite{troelstra:98,kohlenbach:08} for more about
realizability, and \cite{friedman:78} for the A-translation.

G\"odel's Dialectica interpretation provides an alternative route to
this result. In fact, one obtains a stronger conclusion, namely
that the correctness of the witnessing term can be proved in a
quantifier-free fragment $\na{PR^\omega}$ of $\na{HA^\omega}$. To
each formula $\ph$ in the language of arithmetic, one inductively
assigns a formula $\ph^D$ of the form $\ex x \fa y \ph_D(x,y)$, where
$x$ and $y$ are now tuples of variables of appropriate types. Assuming
$\psi^D = \ex u\fa v \psi_D(u,v)$, the assignment is defined as follows:
\begin{alignat*}{2}
&\theta^D & \; & \equiv \; \theta, \; \mbox{if $\theta$ is atomic} \\
& (\varphi \land \psi)^D & \; & \equiv \; \ex{x,u} \fa{y, v}
(\varphi_D\land \psi_D) \\
& (\varphi \lor \psi)^D & \; & \equiv \;  \ex{z} \fa{y, v} 
 (\isleft{z} \land \varphi_D(\eltl{z},y) \mathop{\lor} \\
& & \; & \quad \quad \quad (\isright{z} \land
 \psi_D(\eltr{z}, v))) \\
& (\varphi \limplies \psi)^D & \; & \equiv \; \ex{U,Y} \fa {x,v}
(\varphi_D(x,Y(x,v)) \limplies \psi_D(U(x),v)) \\
& (\fa z \varphi(z))^D & \; & \equiv \; \ex X \fa {z, y}
\varphi_D(X(z),y,z) \\
& (\ex z \varphi(z))^D & \; & \equiv \; \ex{z,x} \fa{y}
\varphi_D(x,y,z)
\end{alignat*}
The clause for implication is the most interesting among these, and
can be understood as follows: from a witness, $x$, to the hypothesis,
$U(x)$ is supposed to return a witness to the conclusion; and given a
purported counterexample, $v$, to the conclusion, $Y(x,v)$ is supposed
to return a counterexample to the hypothesis. Since we have defined
$\lnot \varphi$ to be $\varphi \limplies \bot$, notice that $(\lnot
\ph)^D$ is $\ex Y\fa x\lnot \varphi_D(x,Y(x))$.

The Dialectica interpretation of $\fa x \lnot \lnot \ex y
R(x,y)$ is $\ex Y \fa x \lnot \lnot R(x,Y(x))$, which is
intuitionistically equivalent to $\ex Y \fa x R(x,Y(x))$, given the
decidability of primitive recursive relations. One can show that from
a proof of $\ph$ in $\na{HA}$, one can extract a term $F$ such that
for every $x$, $\na{PR^\omega}$ proves $\ph_D(x,F(x))$, once again
yielding Theorem~\ref{main:witnessing:thm}.

\section{Some double-negation translations}
\label{double:negation:section}

We have seen that one can use modified realizability or the Dialectica
interpretation to extract an algorithm from a proof of a $\Pi_2$
statement in classical arithmetic, modulo a method of reducing
classical arithmetic to intuitionistic arithmetic. Double negation
translations provide the latter.

A formula is said to be {\em negative} if it does not involve
$\exists$ or $\lor$ and each atomic formula $A$ occurs in the form
$\lnot A$; in other words, the formula is built up from negated atomic
formula using $\forall$, $\land$, $\limplies$, and $\bot$. Over
minimal logic, negative formulas are stable under double negation,
which is to say, if $\ph$ is any negative formula, then $\na{HA}$
proves that $\lnot \lnot \ph$ is equivalent to $\ph$ (see, for
example, \cite{troelstra:schwichtenberg:00}). 

The G\"odel-Gentzen double-negation translation maps an arbitrary
first-order formula $\ph$ to a negative formula, $\ph^N$:
\begin{alignat*}{2}
& \bot^N & \; & \equiv \; \bot\\
& \theta^N & \; & \equiv \; \lnot\lnot \theta, \;\mbox{if $\theta$ is atomic}\\
& (\varphi \land \psi)^N & \; & \equiv \; \varphi^N \land \psi^N\\
& (\varphi \lor \psi)^N & \; & \equiv \; \lnot(\lnot \varphi^N \land \lnot
\psi^N)\\
& (\ph \limplies \psi)^N & \; & \equiv \; \ph^N \limplies \psi^N\\
& (\fa x \varphi)^N & \; & \equiv \; \fa x \varphi^N\\
& (\ex x \varphi)^N & \; & \equiv \; \lnot \fa x \lnot \varphi^N
\end{alignat*}
The translation has the following properties:
\begin{theorem}
For any formula $\ph$ and set of sentences $\Gamma$:
\begin{enumerate}
\item Classical logic proves $\ph \liff \ph^N$
\item If $\ph$ is provable from $\Gamma$ in classical logic, then
  $\ph^N$ is provable from $\Gamma^N$ in minimal logic. 
\end{enumerate}
\end{theorem}
Since the $\na{HA}$ proves the $\cdot^N$ translations of its own
axioms, we have as a corollary:
\begin{corollary}
If $\na{PA}$ proves $\ph$, then $\na{HA}$ proves $\ph^N$.
\end{corollary}
In fact, one can strengthen the corollary in three ways:
\begin{enumerate}
\item Since $\na{HA}$ proves $\lnot \lnot \theta \limplies \theta$ for
  atomic formulas $\theta$, one can define $\theta^N$ to be $\theta$.
\item Assuming the language of $\na{HA}$ includes, say, symbols
  denoting the primitive recursive functions, every negated atomic
  formula, $\lnot \theta$, has an atomic equivalent, $\bar \theta$; so
  one can define $(\lnot \theta)^N$ to be $\bar \theta$.
\item The theorem remains true if one replaces $\na{HA}$ by a suitable
  variant, $\na{HA'}$, based on minimal logic.  
\end{enumerate}
These considerations hold in the theorems that follow as well.

The reason to be concerned about negations is that they are
undesirable with respect to the two computational interpretations
given in the last section, since they lead to the use of more
complicated types in the resulting terms of $\PR$. There is a variant
of the double-negation translation known as the \emph{Kuroda
  translation} that fares slightly better in this regard: for any
formula $\ph$, let $\ph^{\mathit{Ku}}$ denote the result of
doubly-negating atomic formulas, and adding a double negation
\emph{after} each universal quantifier, and, finally, adding a
double-negation to the front of the formula. Then we have:
\begin{theorem}
  For every formula $\ph$, $\ph^{\mathit{Ku}} \liff \ph^N$ is provable
  in minimal logic. Hence $\na{PA}$ proves $\ph$ if and only if
  $\na{HA}$ proves $\ph^{\mathit{Ku}}$.
\end{theorem}
\noindent Note that intuitionistic logic, rather than minimal logic,
is required in the conclusion.

Late in 2005, Grisha asked whether a version of the Dialectica
interpretation designed by Shoenfield \cite{shoenfield:01}, for
classical arithmetic, could be understood as a composition of the
usual Dialectica interpretation together with a double-negation
translation. I set the question aside and solved it a few months later
\cite{avigad:unp:k}, only to find that Ulrich Kohlenbach and Thomas
Streicher had solved it more quickly \cite{streicher:kohlenbach:07}.
In a way that can be made precise, the Shoenfield translation
corresponds to the following version of the double-negation
interpretation (itself a variant of a translation due to Krivine),
expressed for a basis involving the connectives $\lnot$, $\land$,
$\lor$, and $\forall$. We define $\ph^{\mathit{Kr}}$ to be $\lnot
\ph_{\mathit{Kr}}$, where $\ph_{\mathit{Kr}}$ is defined recursively
by clauses below.  It helps to keep in mind that $\ph_{\mathit{Kr}}$
is supposed to represent the \emph{negation} of $\ph$:
\begin{alignat*}{2}
  & \theta_{\mathit{Kr}} & \; & \equiv \; \lnot \theta, \; 
    \mbox{if $\theta$ is atomic} \\
  & (\lnot \ph)_{\mathit{Kr}} & \; & \equiv \; \lnot \ph_{\mathit{Kr}} \\
  & (\ph \land \psi)_{\mathit{Kr}} & \; & \equiv \; \ph_{\mathit{Kr}} \lor \psi_{\mathit{Kr}} \\
  & (\ph \lor \psi)_{\mathit{Kr}} & \; & \equiv \; \ph_{\mathit{Kr}} \land \psi_{\mathit{Kr}} \\
  & (\fa x \ph)_{\mathit{Kr}} & \; & \equiv \; \ex x \ph_{\mathit{Kr}}
\end{alignat*}
Note that we can eliminate either $\vee$ or $\wedge$ and retain a
complete set of connectives, but including them both is more
efficient. Formulas of the form $\ex x \ph$, however, have to be
expressed as $\lnot \fa x \lnot \ph$ to apply the translation.
\begin{theorem}
For every formula $\ph$, $\ph^{\mathit{Kr}} \liff \ph^N$ is provable in minimal
logic. Hence $\na{PA}$ proves $\ph$ if and only if $\na{HA'}$
proves $\ph^{\mathit{Kr}}$.
\end{theorem}

The $\cdot^{\mathit{Kr}}$-translation is particularly good when it
comes to formulas in negation-normal form; it only adds two
quantifiers for each existential quantifier, as well as one at the
beginning. But one can do even better \cite{avigad:00}. Taking
advantage of the classical negation operator, now $\ph^M$ is defined
to be $\lnot (\inot \ph)_M$, where the map $\psi \mapsto \psi_M$ is
defined recursively as follows:
\begin{alignat*}{2}
& \theta_M & \; & \equiv \; \theta, \; \mbox{if $\theta$ is atomic or
  negated atomic} \\
& (\ph \lor \psi)_M & \; & \equiv \; \ph_M \lor \psi_M \\ 
& (\ph \land \psi)_M & \; & \equiv \; \ph_M \land \psi_M \\
& (\ex x \ph)_M & \; & \equiv \; \ex x \ph_M \\
& (\fa x \ph)_M & \; & \equiv \; \lnot \ex x (\inot \ph)_M.
\end{alignat*}
Once again, we have
\begin{theorem}
  For every formula $\ph$ in negation-normal form, $\ph^M \liff \ph^N$
  is provable in minimal logic. Hence $\na{PA}$ proves $\ph$ if and
  only if $\na{HA'}$ proves $\ph^M$.
\end{theorem}
The $\cdot^M$-translation is extremely efficient with respect to
negations, introducing, roughly, one at the beginning of the formula,
and one for every quantifier alternation after an initial block of
universal quantifiers. Alternatively, can define $(\ph \land \psi)_M$
in analogy to $(\fa x \ph)_M$, as $\lnot ((\inot \ph)_M \lor (\inot
\psi)_M)$.  This gives the translation the nice property that given
the formulas $\ph^M$ and $(\inot \ph)^M$, one is the negation of the
other. But there is a lot to be said for keeping negations to a
minimum.

Of course, the $\cdot^M$-translation extends to all classical formulas
by identifying them with their canonical negation-normal form
equivalents.  Since the translation relies on the negation-normal form
representation of classical formulas, it shares many nice properties
with a more complicated double-negation translation due to
Girard~\cite{girard:91}. It is this translation that I will use, in
the next section, to provide an efficient computational interpretation
of classical arithmetic.

Let me close with one more translation, found in \cite{avigad:01c},
which is interesting in its own right. For reasons that will become
clear later on, I will call it ``the awkward translation'': if $\ph$
is any formula in negation-normal form, let $\ph^\awk$ denote $\lnot
(\inot \ph)$.

\begin{theorem}
\label{awk:prop}
  For any formula $\ph$, $\ph^N
  \limplies \ph^\awk$ is provable in minimal logic. Hence,
  $\na{PA}$ proves $\ph$ if and only if $\na{HA'}$ proves
  $\ph^\awk$.
\end{theorem}

\begin{proof}
  Once can show by induction that if $\psi$ is any formula in
  negation-normal form, then $\psi \limplies \psi^N$ is provable in
  minimal logic. So minimal logic proves that $\inot \ph$
implies $(\inot \ph)^N$, and hence that $\lnot (\inot \ph)^N$ implies
$\ph^\awk$.  But since $\inot \ph$ is classically
equivalent to $\lnot \ph$, $\lnot (\inot \ph)^N$ is equivalent to
$\lnot \lnot \ph^N$, which is implied by $\ph^N$. 
\end{proof}

The $\cdot^\awk$-translation is almost absurdly efficient with respect
to negations: the one classical negation on the inside adds no
negations at all (recall that in arithmetic, negated atomic formulas
have atomic equivalents), and the translation adds only one negation
on the outside. But the attentive reader will have noticed that the
first assertion in Theorem~\ref{awk:prop} is slightly weaker than the
corresponding assertions in the the theorems that precede it: only one
direction of the equivalence is minimally valid. We will see, in
Section~\ref{herbrand:section}, that this means that the translation
fares very poorly with respect to modus ponens, making it impossible
to translate ordinary proofs piece by piece.

For pure first-order logic, an alternative proof of
Theorem~\ref{awk:prop} can be found in \cite{avigad:01c}. Benno van
den Berg (personal communication) later hit upon this same
translation, independently. In \cite{avigad:01c}, I claimed that with
intuitionistic logic in place of minimal logic, the result is a
consequence of a characterization of Glivenko formulas due to Orevkov,
described in a very nice survey \cite[Section 3.2.5]{mints:91} of
Russian proof theory by Grisha.\footnote{There are typographical
  errors on page 401 of that paper, which Grisha has corrected for me.
  The last class eight lines from the bottom of the page should be
  $\{\rightarrow^-,\inot^-,\vee+,\exists^+\}$; the last class seven
  lines from the bottom of the page should be
  $\{\rightarrow^-,\inot^-,\vee+,\rightarrow^+,\forall^+\}$; and the
  first class at the bottom line should be $\{\rightarrow^+,
  \inot^+,\vee^-\}$.} That seems to be incorrect; but van den Berg and
Streicher have pointed out to me that in that case the result
follows a theorem due to Mints and Orevkov \cite[page 404, paragraph
4]{mints:91}.
 
\section{Interpreting classical arithmetic}
\label{realizability:section}

We now obtain direct computational interpretations of classical
arithmetic simply by combining the $\cdot^M$-translation of
Section~\ref{double:negation:section} with the computational
interpretations of $\na{HA'}$ given in
Section~\ref{preliminaries:section}. One annoying consequence of the
use of the classical negation operator in the $\cdot^M$-translation is
that it is impossible to carry out the translation of a formula $\ph$
from the inside out: depending on the context in which a subformula
$\psi$ occurs, the computational interpretation of the full formula
may depend on either the computational interpretation of $\psi$ or the
computational interpretation of $\inot \psi$. In practice, then, it is
often more convenient to carry out the interpretation in two steps,
applying the $\cdot^M$-translation first, and then one of the two
computational interpretations described in
Section~\ref{preliminaries:section}. Nonetheless, it is interesting to
see what happens when the steps are composed, which is what I will do
here. Both translations apply to formulas in negation-normal form, and
we can assume that negated atomic subformulas are replaced by their
atomic equivalents.

As in Section~\ref{preliminaries:section}, the appropriate version of
classical realizability is defined relative to a fixed primitive
recursive predicate $A(y)$. Most of the clauses look just like
ordinary modified realizability:
\begin{alignat*}{2}
  & a \realizes \theta & \; & \equiv \; \theta, \; \mbox{if $\theta$ is atomic} \\
  & a \realizes \ph \land \psi & \; & \equiv \; ((a)_0 \realizes \ph) \land
  ((a)_1 \realizes \ph) \\
  & a \realizes \ph \lor \psi & \; & \equiv \; ((\isleft{a} \land \eltl{a}
  \realizes
  \ph) \mathop{\lor} \\
  & & \; & \quad \quad \quad (\isright{a} \land \eltr{a} \realizes \psi)) \\
  & a \realizes \ex x \ph(x) & \; & \equiv \; (a)_1 \realizes
  \ph((a)_0) \\
\intertext{The only slightly more complicated clause is the one for
  the 
universal quantifier. Take $a \; \mathop{\mathit{refutes}} \; \ph$ to
be the formula $\fa b (b \realizes \ph \limplies A(a(b)))$.}
  & a \realizes \fa x \ph(x) &\; & \equiv \; a \;
  \mathop{\mathit{refutes}} \; \ex x
  \inot \ph(x)
\end{alignat*}
One can then straightforwardly extract, from any proof of a formula
$\ph$ in classical arithmetic, a term $a$ that refutes $\inot \ph$.
But now notice that the identity function realizes $\fa x \bar A(x)$,
where $\bar A$ is the negation of $A$; so, from a proof of $\ex y
A(y)$ in classical arithmetic, since the identity function realizes
$\fa x \bar A(x)$, one obtains a term $a$ satisfying $A(a)$. This
provides a direct proof of Theorem~\ref{main:witnessing:thm}. Details
can be found in \cite{avigad:00}. A more elaborate realizability
relation, based on the A-translation, can be found in
\cite{berger:et:al:02}.

The corresponding variant of the Dialectica translation is similarly
straightforward. As with the Shoenfield variant
\cite{shoenfield:01,avigad:unp:k}, each formula $\ph$ is mapped to a
formula $\ph^{D'}$ of the form $\fa x \ex y \ph_{D'}(x,y)$, where $x$
and $y$ are sequences of variables.  Assuming $\psi^{D'}$ is $\fa u
\ex v \psi_{D'}(u,v)$, the translation is defined recursively, as
follows:
\begin{alignat*}{2}
  & \theta^{D'} & \; & \equiv \; \theta, \; \mbox{if $\theta$ is atomic} \\
  & (\ph \land \psi)^{D'} & \; & \equiv \; \fa{x,u} \ex{y,v}
  (\ph_{D'}(x,y) \land \psi_{D'}(u,v)) \\
  & (\ph \lor \psi)^{D'} & \; & \equiv \; \fa {x,u} \ex {y,v}
  (\ph_{D'}(x,y) \lor
  \psi_{D'}(u,v)) \\
  & (\fa z \ph)^{D'} & \; & \equiv \; \fa {z,x} \ex y \ph_{D'}(x,y)
  \intertext{This time, it is the clause for the existential
    quantifier that is slightly more complicated. If $(\inot
    \ph(z))^{D'}$ is $\fa r \ex s (\inot \ph)_{D'}(z,r,s)$, define}
  & (\ex z \ph)^{D'} & \; & \equiv \; \fa S \ex {z,r} \lnot (\inot
  \ph)_{D'} (z,r,S(z,r)).
\end{alignat*}
On can interpret this as saying that for any function $S(z,r)$ that
purports to witness $\fa {z,r} \ex s (\inot \ph)_{D'}(z,r,s)$, there
are a $z$ and an $r$ denying that claim. Once again, one can
straightforwardly extract, from any proof of a formula $\ph$ in
classical arithmetic, a term $a$ satisfying $\fa x \ph_{D'} (x,a(x))$.
This provides another direct proof of
Theorem~\ref{main:witnessing:thm}.

\section{Back to the awkward translation} 
\label{herbrand:section}

I would now like to come back to the ``awkward translation,''
discussed at the end of Section~\ref{double:negation:section}. I will
do this via what at first might seem to be a digression through
Kreisel's \emph{no-counterexample interpretation}. Let $\ph$ be a
formula in prenex form, for example,
\[
\ex x \fa y \ex z \fa w \theta(x,y,z,w).
\]
The \emph{Herbrand normal form} $\ph^H$ of $\ph$ is obtained by
replacing the universally quantified variables of $\ph$ by function
symbols that depend on the preceding existential variables, to obtain
\[
\ex {x,z} \theta(x,f(x),z,g(x,z)).
\]
It is not hard to check that, in classical logic, $\ph$ implies
$\ph^H$. Thus, by a slight variant of
Theorem~\ref{main:witnessing:thm} (relativizing it to function
symbols), if classical arithmetic proves $\ph$, there will be terms
$F_1(f,g)$ and $F_2(f,g)$ of $\PR$ such that $\na{HA^\omega}$ proves
\[
\fa {f, g} \theta(F_1(f,g),f(F_1(f,g)),F_2(f,g),g(F_1(f,g),F_2(f,g))).
\]
Think of $f$ and $g$ as providing purported counterexamples to the
truth of $\ph$, so that $F_1$ and $F_2$ effectively foil such
counterexamples. The no-counterexample interpretation is simply the
generalization of this transformation to arbitrary prenex formulas.

One need not invoke Herbrand normal form to arrive at the previous
conclusion.  One can check that if $\ph$ is a prenex formula, the
no-counterexample interpretation of $\ph$ is essentially just the
Dialectica interpretation of $\ph^\awk$, so the result follows from
Theorem~\ref{awk:prop} as well.

The no-counterexample interpretation can be viewed as a computational
interpretation of arithmetic. But, in a remarkable article
\cite{kohlenbach:99}, Kohlenbach has shown that it is not a very
\emph{modular} computational interpretation, in the sense that it does
not have nice behavior with respect to modus ponens. To make this
claim precise, note that the set of (terms denoting) primitive
recursive functionals, $\PR$, can be stratified into increasing
subsets $\PR_n$, in such a way that any finite fragment of $\na{HA}$
has a Dialectica interpretation (or modified realizability
interpretation) using only terms in that set. Kohlenbach
\cite[Proposition 2.2]{kohlenbach:99} shows:
\begin{theorem}
\label{kohlenbach:theorem}
For every $n$ there are sentences $\ph$ and $\psi$ of arithmetic such
that:
\begin{enumerate}
\item $\ph$ is prenex.
\item $\psi$ is a $\Pi_2$ sentence, that is, of the form $\fa x \ex y
  R(x,y)$ for some primitive recursive relation $R$.
\item Primitive recursive arithmetic proves $\ph$.
\item $\na{PA}$ proves $\ph \limplies \psi$.
\item $\ph$ and every prenexation of $\ph \limplies \psi$ has a
  no-counterexample interpretation with functionals in $\PR_0$.
\end{enumerate}
But:
\begin{enumerate}
\item[6.] There is no term $F$ of $\PR_n$ which satisfies the
  no-counterexample interpretation of $\psi$; that is, there is no
  term $F$ such that $\fa x R(x,F(x))$ is true in the standard model
  of arithmetic.
\end{enumerate}
\end{theorem}
\noindent Theorem~\ref{kohlenbach:theorem} shows that there is no
straightforward way to combine witnesses to the no-counterexample
interpretations of $\ph$ and $\ph \limplies \psi$, respectively, to
obtain a witness to the no-counterexample interpretation of $\psi$.

The problem with the awkward translation is that, similarly, it may
behave poorly with respect to modus ponens. Consider pure first-order
logic with a single predicate symbol, $A(x)$. Then there are formulas
$\ph$ and $\psi$ such that $\psi^\awk$ doesn't follow from $\ph^\awk$
and $(\ph \limplies \psi)^\awk$ in minimal logic: just take $\ph$ to
be the formula $\fa x A(x)$ and $\psi$ to be $\bot$. In that case,
$\ph^\awk \land (\ph \limplies \psi)^\awk \limplies \psi^\awk$ is
equivalent, over minimal logic, to the double-negation shift, $\fa x
\lnot \lnot A(x) \limplies \lnot \lnot A(x)$, which is not even
provable in intuitionistic logic.

When I showed the awkward translation to Grisha, he remarked right
away that its behavior has something to do with Kohlenbach's result.
At the time, I had no idea what he meant; but writing this paper
finally prodded me to sort it out. Grisha was right:
Theorem~\ref{kohlenbach:theorem} can, in fact, be used to show that
the awkward translation does not provide a modular translation of
Peano arithmetic to Heyting arithmetic, in the following sense.

\begin{theorem}
\label{awk:theorem}
For any fragment $T$ of $\na{HA}$, there are formulas $\ph$ and $\psi$
such that the following hold:
\begin{enumerate}
\item $\na{PA}$ proves $\ph$ and $\ph \limplies \psi$, but
\item $\na{T}$ together with $\ph^\awk$ and $(\ph \limplies
  \psi)^\awk$ does not prove $\psi^\awk$.
\end{enumerate}
\end{theorem}
\noindent This shows that modus ponens fails under the
$\cdot^\awk$-translation, in a strong way. Since $\na{PA}$ proves
$\ph$ and $\ph \limplies \psi$, it also proves $\psi$, and so by
Theorem~\ref{awk:prop}, $\na{HA}$ proves $\psi^\awk$. But having the
translation $\ph^\awk$ and $(\ph \limplies \psi)^\awk$ may not help
much in obtaining such a proof of $\psi^\awk$; indeed, obtaining a
proof of $\psi^\awk$ from $\ph^\awk$ and $(\ph \limplies
\psi)^\awk$ may be no easier than simply proving $\psi^\awk$ outright.

\begin{proof}
  Given $T$, first let $n$ be large enough so that the Dialectica
  interpretation of $T$ uses only terms in $\PR_n$, and then let $\ph$
  and $\psi$ be as in Theorem~\ref{kohlenbach:theorem}. If there were
  a proof of $\psi^\awk$ from $\ph^\awk$ and $(\ph \limplies
  \psi)^\awk$ in $T$, applying the Dialectica interpretation, one
  would obtain terms witnessing the Dialectica interpretation of
  $\psi^\awk$ from terms witnessing the Dialectica interpretations of 
  $\ph^\awk$ and $(\ph \limplies \psi)^\awk$. But the Dialectica
  interpretation of $\ph^\awk$ is the no-counterexample interpretation
  of $\ph$, and it is not hard to check that the Dialectica
  interpretation of $(\ph \limplies \psi)^\awk$ is the
  no-counterexample interpretation of one of the prenexations of $\ph
  \limplies \psi$. Thus there would be a witness to the
  no-counterexample interpretation of $\psi$ in $\PR_n$, contrary to
  the choice of $\ph$ and $\psi$.
\end{proof}

\section{Conclusions}
\label{conclusions:section}

As noted in the introduction, conventional wisdom holds that classical
logic is ``nondeterministic,'' in that different ways of extracting
algorithms from classical proofs can yield different results.
Sometimes, however, nondeterminacy is unavoidable. For example, even
minimal logic can prove $\ex x A(x) \land \ex y A(y) \limplies \ex z
A(z)$, and any computational interpretation of this formula will have
to choose either $x$ or $y$ to witness the conclusion. The difference
is that this sentence is typically not taken to be an \emph{axiom} of
minimal logic; rather, there are two axioms, $\ph \land \psi \limplies
\ph$ and $\ph \land \psi \limplies \psi$, and any proof of the
sentence has to choose one or the other. In contrast, standard calculi
for classical logic provide cases where there are multiple choices of
witnesses, with no principled reason to choose one over the other. For
example, starting from canonical proofs of $A(a) \limplies \ex x A(x)$
and $A(b) \limplies \ex x A(x)$, one can weaken the conclusions to
obtain proofs of $\ph \limplies (A(a) \limplies \ex x A(x))$ and
$\lnot \ph \limplies (A(b) \limplies \ex x A(x))$, for an arbitrary
formula, $\ph$. Then, using the law of the excluded middle, $\ph \lor
\lnot \ph$, one can combine these to obtain a proof of $A(a) \land
A(b) \limplies \ex x A(x)$ where there is little reason to favor $a$
or $b$ as the implicit witness to the existential quantifier. (This
example is essentially that given by Lafont \cite[p.
150]{girard:89}.)

This shows that standard classical calculi are nondeterministic in a
way that proofs in intuitionistic and minimal logic are not. Girard
\cite{girard:91} has neatly diagnosed the source of the
nondeterminacy, and has provided a calculus for classical logic that
eliminates it by forcing the prover to make an explicit choice in
exactly those situations where an ambiguity would otherwise
arise.\footnote{One way to understand what is going on is to notice
  that in minimal logic there are two distinct ways of proving $\lnot
  \lnot (\ph \land \psi)$ from $\lnot \lnot \ph$ and $\lnot \lnot
  \psi$, and this inference is needed to verify the classical axiom
  $\lnot \lnot \theta \limplies \theta$ under the double-negation
  translation.  Another way is to notice that since, in classical
  logic, $\ph$ and $\lnot \lnot \ph$ are equivalent, there is little
  reason to favor $\ph$ or $\lnot \ph$ in situations where
  intuitionistic logic treats them differently. I am grateful to
  Thomas Streicher for these insights.} But note that the
realizability interpretation of Section~\ref{realizability:section}
also avoids this nondeterminacy; the translation procedure described
in \cite{avigad:00} is fully explicit and unambiguous. With respect to
the example discussed in the last paragraph, the interpretation
chooses $a$ or $b$ based on the logical form of $\ph$.\footnote{More
  generally, what breaks the symmetry alluded to in the previous
  footnote is that $\ph$ and $\lnot \lnot \ph$ have the same
  negation-normal form, which is distinct from (and dual to) that of
  $\lnot \ph$.} Indeed, the results of \cite{avigad:00} show that the
witnesses obtained in this way coincide with those obtained using a
natural class of cut-elimination procedures.

Yet another response to the example above is that of Urban and Bierman
\cite{urban:00,urban:bierman:01}, who simply embrace the
nondeterminism as an inherent part of the computational interpretation
of classical logic. In fact, Urban \cite{urban:00} has provided a
nondeterministic programming language to interpret classical logic in
a natural way. Passing through a double-negation interpretation, as we
have done here, amounts to making specific choices to resolve the
nondeterminism. It is an open-ended conceptual problem to understand
which deterministic instances of the general nondeterministic
algorithms can be realized in such a way.

At this point, however, we should be clearer as to the goals of our
analysis. Ordinary mathematical proofs are not written in formal
languages, and so the process of extracting an algorithm from even a
rather constructive mathematical argument can involve nondeterminism
of sorts. And, despite some interesting explorations in this direction
\cite{berger:et:al:01}, it is far from clear that classical arithmetic
can be used as an effective programming language in its own right. But
formal methods are actively being developed in support of software
verification \cite{hoare:03}, and a better understanding of the
computational content of classical logic may support the development
of better logical frameworks for that purpose
\cite{pfenning:01}. Formal translations like the ones described here
have also been effective in ``proof mining,'' the practice of using
logical methods to extract mathematically useful information from
nonconstructive proofs
\cite{avigad:et:al:unp,kohlenbach:08,kohlenbach:oliva:03a}.

Grisha's work has, primarily, addressed the general foundational
question as to the computational content of classical methods. In that
respect, the general metatheorems described here provide a satisfying
answer: for the most part, classical mathematical reasoning
\emph{does} have computational content, which is to say, algorithms
can be extracted from classical proofs; but by suppressing
computational detail, the proofs often leave algorithmic detail
underspecified, rendering them amenable to different
implementations. Grisha's work has thus contributed to an
understanding of the computational content of classical arithmetic
that is mathematically and philosophically satisfying, providing a
solid basis for further scientific research.


\newcommand{\nameindex}[1]{}


\end{document}